\documentstyle[11pt]{article}    

\author {BRUNO FABRE}   

\begin{document}  
\def\B{ \mbox{I\hspace{-.15em}B}}  
\def\dem{ {\it D\'emonstration}.\\ } 
\def\CQFD{ \hbox{\vrule height 5pt depth 5pt width 5pt}} 
\def\N{ \mbox{I\hspace{-.15em}N}} 
\def\Z{ \mbox{Z\hspace{-.3em}Z}} 
\def\R{\mbox{I\hspace{-.17em}R}}   
\def\P{\mbox{I\hspace{-.17em}P}}   
\def\C{ \mbox{l\hspace{-.47em}C}} 
\def\Q{\mbox{l\hspace{-.47em}Q}}   

\newtheorem{thm}{Th\'eor\`eme\hspace{2pt}} 
\newtheorem {prb}{Probl\`eme\hspace{2pt}} 
\newtheorem {pro}{Proposition\hspace{2pt}} 
\newtheorem {conj}{Conjecture\hspace{2pt}}  
\newtheorem {cor}{Corollaire\hspace{2pt}} 
\newtheorem{defi}{Definition\hspace{2pt}} 
\newtheorem{asum}{Hypoth\`ese\hspace{2pt}} 
\newtheorem {rem}{Remarque\hspace{2pt}} 
\newtheorem {lem}{Lemme\hspace{2pt}}  
\newtheorem {exm}{Exemple\hspace{2pt}}    

\baselineskip=12.72pt   

 {\center\LARGE\bf{Sur la transformation d'Abel-Radon \\ {des courants localement 
r\'esiduels}}}  
\\
\\{BRUNO FABRE}
\\

 {\bf Abstract.} 
We give in this note a generalisation of the following theorem of Henkin
and Passare (cf. \cite{Princeton} and \cite{HenPas}) : Let be $Y$ an analytic
subvariety of pure codimension $p$ in a linearly $p-$concave domain $U$, and
$\omega$ a meromorphic $q-$form ($q>0$) on $Y$; if the Abel-Radon transform ${\cal
R}(\omega\wedge[Y])$, which is meromorphic on $U^*$, has a meromorphic prolongation
to ${\tilde U}^*$ containing $U^*$, then $Y$ extends as an analytic subvariety $\tilde Y$ of $\tilde U$, and
$\omega$ as a meromorphic form on $\tilde Y$. We show the analogous
statement when we replace $\omega\wedge[Y]$ by a current $\alpha$ of a more general
type, called {\it locally residual}, if $\alpha$ is of bidegree
$(N,1)$, or $(q+1,1),0<q<N$ in the particular case where ${\cal R}(\alpha)=0$. 


{\bf Introduction.} 

 On commence dans les deux premiers paragraphes par introduire les deux objets que nous utiliserons par la suite,
 dans un contexte assez g\'en\'eral: premi\`erement, les courants localement r\'esiduels, et deuxi\`ement, la
transformation 
d'Abel-Radon, dont la transformation d'Abel (cf.\cite{Griffiths} et \cite{Fabre2}) est un cas 
particulier. On montre pour la transformation d'Abel-Radon certaines propri\'et\'es 
de stabilit\'e, en particulier que la transformation d'un courant localement 
r\'esiduel est encore un courant localement r\'esiduel. 

 Le probl\`eme principal \'etudi\'e dans cette note peut alors se formuler de la mani\`ere suivante. On se donne un
ouvert de l'espace projectif complexe $\P_N$, lin\'eairement $p-$concave, i.e. r\'eunion de $p-$plans complexes.
On peut donc d\'efinir la transformation d'Abel-Radon ${\cal R}$, qui associe \`a un courant de bidegr\'e
$(r,s)$ sur $U$ un courant de bidegr\'e $(r-p,s-p)$ sur $U^*$, l'ouvert de la grassmanienne $G(p,N)$
correspondant aux $p-$plans contenus dans $U$. Si $\alpha$ est un courant localement r\'esiduel de bidegr\'e
$(N,p)$ sur $U$, ${\cal R}(\alpha)$ est une $n-$forme m\'eromorphe sur $U^*$ ($n=N-p$). Alors, si ${\cal R}(\alpha)$
se prolonge sur un domaine ${\tilde U}^*$ contenant $U^*$, peut-on prolonger $\alpha$ dans $\tilde U$ comme courant
localement r\'esiduel ?

Nous donnons la r\'eponse positive pour $p=1$. On commence par
\'etudier, \'etant donn\'ee la projection standard $\phi:D\times\C_y\to D$, o\`u $D$ est un domaine de $\P_n$, la
trace $\phi_*(\alpha)$ des courants localements r\'esiduels $\alpha$ \`a support $\phi-$propre. On montre que
$\alpha$ peut se reconstruire, comme courant {\it r\'esiduel}, \`a partir de ses traces $u_k:=\phi_*(\alpha
y^k),k\in \N$; puis que si les $u_k$ se prolongent m\'eromorphiquement sur un domaine plus grand $\tilde D$,
$\alpha$ se prolonge sur $\tilde D\times\C$ comme courant r\'esiduel. 

  On applique ensuite ce th\'eor\`eme pour d\'emontrer le th\'eor\`eme suivant:
 
Si $U$ est un domaine de $\P_N$ r\'eunion de droites complexes, et $\alpha$ un courant localement r\'esiduel de
bidegr\'e $(N,1)$ sur $U$, si la transformation d'Abel-Radon de $\alpha$, m\'eromorphe sur $U^*$, se prolonge
m\'eromorphiquement sur un domaine ${\tilde U}^*$ contenant $U^*$, alors $\alpha$ se prolonge comme courant localement
r\'esiduel sur le domaine $\tilde U:=\cup_{t\in {\tilde U}^*}{\Delta_t}$ contenant $U$.

  Ce th\'eor\`eme est une g\'en\'eralisation du th\'eor\`eme d'Henkin-Passare 
(cf. \cite{HenPas} et \cite{Princeton}), qui correspond au cas o\`u $\alpha$ est de la forme $\omega\wedge[Y]$, au
cas o\`u $\alpha$ est \`a "support non-r\'eduit", et r\'esoud pour $N=2$ un probl\`eme formul\'e dans
\cite{NouHen}. Ce th\'eor\`eme d'Henkin-Passare est lui-m\^eme une g\'en\'eralisation du th\'eor\`eme de Griffiths,
qui correspond au cas o\`u la transformation d'Abel-Radon est nulle. 
  Dans le cas o\`u la transformation d'Abel-Radon est nulle,
 on peut d\'eduire du th\'eor\`eme d'injectivit\'e d'Henkin-Gindikin (\cite{GHenkin}) un th\'eor\`eme d'extension
pour tous les bidegr\'es $(q+1,1), 0<q<N$, comme on l'a montr\'e dans \cite{Fabre}, g\'en\'eralisant ainsi le
th\'eor\`eme de Griffiths au cas non r\'eduit. 

On donne une application du th\'eor\`eme d'Henkin-Passare g\'en\'eralis\'e: \'etant donn\'e un
courant localement r\'esiduel de bidegr\'e $(N,1)$ sur un domaine $U\subset \P_N$ lin\'eairement $1$-concave, tel que $U^*$
est connexe, si
$U$ contient un $2-$plan, le courant se prolonge \`a $\P_N$ en un courant localement r\'esiduel. 

On termine par mentionner quelques questions ouvertes naturelles dans ce contexte. En particulier, on montre la difficult\'e
pour g\'en\'eraliser le th\'eor\`eme de Griffiths au cas de courants localement r\'esiduels de bidegr\'es $(q+2,2)$.

{REMERCIEMENTS.} 

Je remercie le professeur Tomassini pour son acceuil 
 \`a l'\'ecole normale sup\'erieure de Pise, o\`u j'ai pu terminer la r\'edaction de 
cet article, ainsi que pour les discussions encourageantes qu'il a bien voulu 
m'accorder. Je remercie d'ailleurs aussi d'autres math\'ematiciens pour des 
conversations fructueuses: G. Henkin, P. Mazet, J.-L. Sauvageot.

\section{Rappels.}

 Il existe une seule mani\`ere d'associer \`a une vari\'et\'e diff\'erentiable orient\'ee $X$ de dimension $n$ et \`a une
$n-$forme continue
\`a support compact $\omega$ sur $X$ un nombre $\int_X{\omega}$,
de sorte que:

i) Pour $X$ fix\'e, l'application $\omega\mapsto \int_X{\omega}$ est lin\'eaire continue;

ii) si $X$ est un point $P$, affect\'e d'un signe $\sigma$, alors $\int_X{\omega}=\sigma\omega(P)$;

iii) Si $Y$ est un ouvert de $X$, \`a bord lisse, de bord
$\partial Y$ compact, orient\'e d'apr\`es l'orientation de $Y$, contenant
le support de $d\omega$, alors
$\int_{Y}{d\omega}=\int_{\partial Y}{\omega}$.

Soit $X$ une vari\'et\'e diff\'erentiable de dimension $n$.
Etant donn\'e un compact $K\subset X$, on note ${\cal D}_{K}$ l'espace des fonctions lisses \`a support
dans
$K$, muni de la topologie d\'efinie par les semi-normes $p_{m, K}:=\sup_{K}{D^m(f)}$ (o\`u $D^m$ est la
d\'eriv\'ee
$m-$i\`eme de $f$),
 avec $0\le m< \infty$. 
On d\'efinit une topologie sur l'espace ${\cal D}$ des fonctions lisses \`a support compact de la mani\`ere
suivante: soit $K_i$ une famille croissante de compacts de r\'eunion $X$. Alors ${\cal D}=\cup_i{{\cal D}_{K_i}}$. On dit
que $U$ est un ouvert de ${\cal D}$ ssi: $U\cap{\cal D}_{K_i}$ est un ouvert de ${\cal D}_{K_i}$ pour tout $i$. 
Pour cette topologie, une suite de
fonctions $f_n$ tend vers
$0$ dans ${\cal D}$ si:

 i) il existe un compact $K$ contenant le support des $f_n$;

ii) $f_n$ tend vers $0$ dans ${\cal D}_K$.

On d\'efinit de mani\`ere analogue une topologie sur l'espace ${\cal D}^p$ des $p-$formes lisses \`a support
compact. 

\begin{defi}
Un courant de degr\'e $p$ est une forme lin\'eaire continue sur ${\cal D}^{n-p}$. Une distribution est un
courant de degr\'e $0$.
\end{defi}

Etant donn\'e un fonction localement int\'egrable $f$, on lui associe une distribution $[f]$ par la
formule: $[f](\phi):=\int_X{f\phi}$.
On peut multiplier un courant $T$ par une forme lisse $\phi$, en posant:
$(\phi\wedge T)(\psi):=T(\phi\wedge\psi)$.
Alors un courant de degr\'e $p$, dans un syst\`eme de coordonn\'ees locales $x_1,\dots,x_n$,  s'\'ecrit
de mani\`ere unique comme $T=\sum_{I=(i_1<\dots<i_p)}{T_Idx^I}$, avec $T_I$ des distributions.

On \'etend les op\'erateurs de d\'erivation (ou champs de vecteurs) aux distributions, de sorte que $D[f]=[Df]$ si $f$ est
une fonction lisse, et $D(fg)=Df[g]+fDg$ si $f$ est une distribution et $g$ une fonction lisse. 
De m\^eme, on
\'etend l'op\'erateur diff\'erentielle
$d$, pour les distributions puis pour les courants : en coordonn\'ees locales, $dT:=\sum_i{\partial_{x_i}Tdx_i}$ et
$d{\sum_I{T_Idx^I}}:=\sum_I{dT_I\wedge dx^I}$.

 Soit $X$ un espace analytique complexe r\'eduit.
 Alors toute
$k-$forme diff\'erentiable se d\'ecompose de mani\`ere unique en $k-$formes de bidegr\'e
$(p,q)$, avec $p+q=k$: $\omega:=\sum_{p+q}{\omega_{p,q}}$; on en d\'eduit que tout courant $T$ de degr\'e
$k$ se d\'ecompose $T=\sum_{p+q=k}{T_{p,q}}$, avec $T_{p,q}$ un courant de bidegr\'e $(p,q)$, i.e. une
$(p,q)-$forme \`a coefficients distributions. L'op\'erateur de diff\'erentiation
$d$ se d\'ecompose alors en somme :
$d=\partial+\overline\partial$, de sorte que $\partial$ associe \`a toute forme (resp. courant) de
bidegr\'e
$(p,q)$ une forme (resp. courant) de bidegr\'e $(p+1,q)$, et
$\overline\partial$ associe \`a toute forme (resp. courant) de bidegr\'e $(p,q)$ une forme (resp.
courant) de bidegr\'e
$(p,q+1)$.

 \section{Courants localements r\'esiduels.} 

\subsection{Valeur principale et r\'esidu.}

Soit $X$ une vari\'et\'e analytique. 

Etant donn\'e une fonction m\'eromorphe $f$ sur $X$,
 on a une distribution $[f]$ bien d\'efinie en dehors de son ensemble
polaire $T$. Si $f$ est localement int\'egrable, cette distribution s'\'etend de mani\`ere canonique
 en distribution sur $X$, comme on l'a vu dans le cas diff\'erentiable. 
Si $f$ n'est pas localement int\'egrable, le fait que $f$ est m\'eromorphe permet de montrer
 l'existence d'un prolongement, comme
courant. Mais on ne sait pas a priori faire un choix "canonique" parmi tous les prolongements possibles. Herrera et
Liebermann ont montr\'e dans
\cite{HerLib} le th\'eor\`eme suivant :
\begin{lem} Soit $g$ une fonction holomorphe qui s'annule sur le p\^ole de $f$.
Alors $\lim_{\epsilon\to 0}{\int_{\vline g\vline=\epsilon}{f\phi}}$ et $\lim_{\epsilon\to 0}{\int_{\vline
g\vline>\epsilon}{f\phi}}$ existent pour $\phi$ \`a support compact, et d\'efinissent des courants, not\'es
respectivement $[f]$ et $Res_T(f)$. De plus, $Res_T(f)=\overline\partial[f]$. 
\end{lem}   

 De 
mani\`ere 
g\'en\'erale, si $\psi$ est une forme lisse, et on notera $[\psi f]:=\psi[f]$. 

Soit $\chi(x):\R^+\to \R^+$ une fonction continue et croissante telle que $\chi(x)=0$ au voisinage de $0$ et $\chi(x)=1$
pour
$x\ge 1$, et soit
$g$ une fonction holomorphe s'annulant sur l'ensemble polaire de $f$.
Alors, on d\'eduit de ce qui pr\'ec\`ede : $\lim_{\epsilon\to 0}{\chi(\vline g\vline/\epsilon)[f]}=[f]$ dans l'espace des
distributions.

\begin{defi} Une distribution $\alpha$ sur $X$ est dite d'extension standard 
si pour tout fonction holomorphe $g$ non nulle sur un domaine $U$ de $X$, et pour toute fonction lisse $\chi(x):\R^+\to
\R^+$ croissante telle que $\chi(x)=0$ au voisinage de $0$ et $\chi(x)=1$ pour
$x\ge 1$, $\lim_{\epsilon\to 0}{\chi(\vline
g\vline/\epsilon)\alpha}=
\alpha$ dans $U$.
\end{defi}

Alors, ce qui pr\'ec\`ede montre que la distribution $[f]$ est d'extension standard.

Donc : 
\begin{lem}
Soit $\alpha$ un courant de bidegr\'e $(q,0)$.
Supposons que :

 i) $\alpha$ est $\overline\partial-$ferm\'e en dehors d'une hypersurface complexe.

ii) $\alpha$ est d'extension standard.

Alors, $\alpha=[\omega]$ pour une certaine $q-$forme m\'eromorphe $\omega$.
\end{lem}

\begin{rem}
L'autre propri\'et\'e caract\'eristique de la valeur principale est la suivante :
si $D$ est une d\'erivation holomorphe, et $f$ une fonction m\'eromorphe, alors $Df$ est encore m\'eromorphe, et
 $[Df]=D[f]$.
\end{rem}

On suppose maintenant $X$ un espace analytique r\'eduit. Herrera et Lieberman (\cite{HerLib}) on montr\'e aussi dans ce cas
l'existence du courant "valeur principale" et du courant r\'esidu associ\'e.

\begin{defi} La $q-$forme m\'eromorphe $\omega$ sur $X$ est dite ab\'elienne si le courant $[\omega]$ sur $Y$ est
$\overline\partial-$ferm\'e.
\end{defi}

Soit $\omega$ une $q-$forme holomorphe, en dehors d'une hypersurface $T$ de $X$, contenant $Sing(X)$.
 On a vu que si $\omega$ est
m\'eromorphe (i.e. admet un prolongement m\'eromorphe \`a $X$), le courant $[\omega]$ se prolonge \`a $X$.
 Henkin et Passare ont
montr\'e dans
\cite{HenPas} la r\'eciproque :

\begin{lem}
\label{HenPas}
Si le courant $[\omega]$ sur $X-T$ se prolonge en un courant $\alpha$ sur $X$, $\omega$ est m\'eromorphe sur $X$.
 De plus, si
$\overline\partial\alpha=0$, alors $\omega$ est ab\'elienne. 
\end{lem} 

\begin{rem} Il existe en g\'en\'eral des courants $\overline\partial-$ferm\'es de bidegr\'e $(q,0)$ \`a support dans $Sing(X)$.
\end{rem}

\subsection{Courants localement r\'esiduels de codimension quelconque.}

Soit $X$ une vari\'et\'e analytique de dimension $n$.

 Supposons donn\'ee $(f_1,\dots,f_p,f_{p+1})$ 
 une suite r\'eguli\`ere de fonctions holomorphes sur $X$ (i.e. $f_1=\dots=f_p=g$ est un sous-ensemble analytique
de $U$ de codimension $p+1$). Alors on peut montrer (cf. \cite{Passare})
que les limites, pour $\phi$ lisse \`a support compact :
$$\lim_{t\to 
0}{\int_{\vline f_1\vline=t\epsilon_1,\dots,\vline 
f_p\vline=t\epsilon_p}{\phi/{(f_1\dots f_p)}}},$$ et
$$\lim_{t\to 
0}{\int_{\vline g\vline\ge t\epsilon_{p+1},\vline f_1\vline=t\epsilon_1,\dots,\vline 
f_p\vline=t\epsilon_p}{\phi/{(f_1\dots f_p g)}}}$$ 
existent et sont ind\'ependantes de $\epsilon:=(\epsilon_1,\dots,\epsilon_{p})\in
\R_{+}^{p}$
(resp. $\epsilon':=(\epsilon_1,\dots,\epsilon_{p+1})\in \R_{+}^{p+1},$) en dehors d'un ensemble de mesure nulle;
et que de plus, ces limites d\'efinissent des courants, not\'es
$\overline\partial[1/f_1]\wedge\dots\wedge\overline\partial[1/f_p]$
et $[1/g]\overline\partial[1/f_1]\wedge\dots\wedge\overline\partial[1/f_p]$,
de bidegr\'es $(0,p)$.
Un tel courant, \'eventuellement multipli\'e par une $q-$forme holomorphe, s'appelera {\it courant r\'esiduel} de
codimension $p$.

\begin{rem}
On sait que sous certaines conditions relatives au "front d'onde", on peut multiplier les distributions. 
Il serait int\'eressant de savoir si, dans le cas o\`u $f_1,\dots,f_p,g$ forment une suite r\'eguli\`erede fonctions
holomorphes, on peut multiplier les distributions
$\partial_{\overline{z_1}}[1/f_1],\dots,\partial_{\overline{z_p}}[1/f_p],[1/g]$. Les notations pr\'ec\'edentes se
trouveraient alors parfaitement justifi\'ees.
\end{rem}

De plus, on a :

i) $\overline\partial([1/g]\overline\partial[1/f_1]\wedge\dots\wedge\overline\partial[1/f_p])=
\overline\partial[1/g]\wedge\overline\partial[1/f_1]\wedge\dots\wedge\overline\partial[1/f_p].$
En particulier, $\overline\partial[1/f_1]\wedge\dots\wedge\overline\partial[1/f_p]$ est $\overline\partial-$ferm\'e.

ii) On a :
$f_i[1/g]\overline\partial[1/f_1]\wedge\dots\wedge\overline\partial[1/f_p]=0$. En particulier, le support est contenu dans
$\{f_1=\dots=f_p=0\}$.

iii) Si le support d'un courant r\'esiduel de codimension $p$ est de codimension $p+1$, alors il est nul.

On en d\'eduit :
\begin{lem} Le support d'un courant
r\'esiduel de codimension $p$ est un ensemble analytique de codimension pure $p$.
\end{lem}

\begin{defi}
Si un 
courant est 
r\'esiduel au voisinage de tout point, il s'appellera {\it localement r\'esiduel} (l.r.).
\end{defi}

Pour r\'esumer :

i) Un courant l.r. de bidegr\'e $(q,0)$ est la
valeur principale associ\'ee \`a une $q-$forme m\'eromorphe.

ii) Etant donn\'e un courant l.r. $\alpha$ $\overline\partial-$ferm\'e de bidegr\'e $(q,p)$, il s'\'ecrit localement sous la
forme
$\omega\wedge\overline\partial[1/f_1]\wedge\dots\wedge\overline\partial[1/f_p],$ avec $\omega$ une $q-$forme holomorphe. On
peut le multiplier par un courant de la forme $[g]$, o\`u $g$ est une forme m\'eromorphe dont le p\^ole coupe le support
de
$\alpha$ proprement. Alors, $\overline\partial([g]\alpha)=\overline\partial[g]\wedge\overline\partial\alpha$ est un courant
l.r. de codimension $p+1$.

Soit $\alpha$ un courant l.r. de bidegr\'e $(q,p)$. 
\begin{defi} L'ensemble polaire de $\alpha$ est le support de $\overline\partial\alpha$.
\end{defi}
En dehors de l'ensemble polaire de $\alpha$ (qui est un ensemble analytique de codimension $p+1$), $\alpha$ s'\'ecrit
localement sous la forme
$\omega\wedge\overline\partial[1/f_1]\wedge\dots\wedge\overline\partial[1/f_n]$, avec $\omega$ une $q-$forme
holomorphe.

 Si $p=n=\dim(X),$, et $\omega$ une $n-$forme holomorphe, $f_1,\dots,f_n$ $n$ fonctions holomorphes se coupant dans un
ensemble {\it fini} de points $(P_i)$. Alors, le courant $\omega\wedge
\overline\partial[1/f_1]\wedge\dots\wedge\overline\partial[1/f_n]$ peut se calculer sur la fonction constante $1$; on
obtient:
$$\omega\wedge
\overline\partial[1/f_1]\wedge\dots\wedge\overline\partial[1/f_n](1)=\sum_{P_i}{ Res_{P_i}{\omega/{f_1\dots 
f_n}}},$$ o\`u
$Res$ est 
le r\'esidu ponctuel de Grothendieck, comme on le voit en \'ecrivant la d\'efinition du courant. 

  D'apr\`es J.-E. Bj\"ork (\cite{Bjork}), on peut caract\'eriser les courants l. r. $\overline\partial-$ferm\'es sur $X$ de
la mani\`ere suivante.
 Soit $\alpha$ un courant de bidegr\'e $(q,p)$.
On suppose : 

1) $\alpha$ \`a support dans un ensemble analytique $Z$ de codimension pure 
$p$;

2) $\overline\partial\alpha=0,$

3) ${\overline{I_Z}}\alpha=0,$ i.e. pour tout germe en $z\in Z$ $f_z$ de fonction holomorphe au voisinage de $z$,
$\overline{f_z}\alpha=0$ au voisinage de $z$;

4) $\alpha$ est d'{\it extension standard}, c'est-\`a-dire que pour toute 
fonction holomorphe $g$ dans un voisinage ouvert $U_z$ d'un point $z\in Z$, telle que 
$g$ ne s'annule sur aucune composante irr\'eductible de $U_z\cap Z$, et tout $\chi:\R^+\to \R^+$ croissante nulle dans un
voisinage de $0$ et \'egale \`a $1$ pour $x$ assez grand, on a
$\lim_{\epsilon\to  0}{\chi(\vline g\vline/\epsilon)\alpha}=\alpha$ sur $U_z$. 

 Alors $\alpha$ est l. r., 
s'\'ecrit donc dans un voisinage $U_z$ d'un 
point arbitraire $z\in U$ sous la forme 
$\omega\wedge {\overline\partial[1/f_1]}\wedge\dots\wedge {\overline\partial[1/f_p]}$, 
avec $\omega$ une $q-$forme holomorphe et
$f_1,\dots,f_p$ des fonctions holomorphes formant une suite r\'eguli\`ere sur $U_z$. 

On voit plus g\'en\'eralement que la caract\'erisation reste valide si l'on remplace la deuxi\`eme condition par :

2') $\alpha$ est $\overline\partial-$ferm\'e en dehors d'un sous-ensemble analytique de codimension $p+1$.

Soit $Y$ un sous-ensemble analytique de $X$. Soit $\omega$ une $q-$forme m\'eromorphe sur $Y$; on a vu ci-dessus qu'on peut
d\'efinir le courant valeur principale $[\omega]$ sur $Y$.
 Alors le courant $\omega\wedge[Y]$, est d\'efini par
$\omega\wedge[Y](\phi):=[\omega](i^*(\phi))$, avec
$i$ l'inclusion
$i:Y\subset X$. 

\begin{lem}
Soit $Y$ un ensemble analytique de dimension pure.
Le courant $\omega\wedge[Y]$ est l. r. .
 \end{lem}

\dem Si $Y$ est de codimension $1$, c'est un r\'esultat classique : $\omega$ est le r\'esidu de Poincar\'e-Leray d'une
$(q+1)-$forme m\'eromorphe $\Psi$ sur $X$ \`a p\^ole logarithmique sur $Y$. Si on \'ecrit $\Psi=\Psi'/g$, avec $g$
coupant $Y$ transversalement et $\Psi'$ holomorphe en dehors de $Y$, on a :
$\omega\wedge [Y]=[1/g]\overline\partial[\Psi']$.
Supposons qu'on ait montr\'e le lemme pour $Y$ de codimension $p-1$. Alors, localement,
$Y$ est contenu dans un ensemble analytique $Y'$ de codimension $p-1$, et le courant $[\omega]$ sur $Y$ d\'efinit un
courant sur $Y'$. On a donc : $\omega\wedge [Y]=[1/g]\overline\partial (\Psi'\wedge[Y']),$ avec $\Psi'$ une certaine forme
m\'eromorphe sur $Y'$. Mais on sait que $\Psi'\wedge[Y']$ est l.r. . $g$ \'etant une fonction m\'eromorphe coupant $Y$
proprement, il en est de m\^eme de $[1/g]\overline\partial (\Psi'\wedge[Y'])=\omega\wedge[Y]$.
\CQFD

\section{Transformation d'Abel-Radon.}

\subsection{Image directe d'un courant.} 

Soit un morphisme analytique $\phi : Y\to X$ entre deux vari\'et\'es analytiques complexes de dimension pure, avec
$\dim(Y)=\dim(X)+p$. On peut associer \`a tout courant $\alpha$ de bidegr\'e 
$(r,s)$ sur $Y$, tel 
que la restriction de $\phi$ au support $Supp(\alpha)$ de $\alpha$ est propre, un 
courant 
$\phi_*(\alpha)$ (l'image directe de $\alpha$ par $\phi$) de m\^eme bidimension sur $X$, donc de bidegr\'e
$(r-p,s-p)$, par la formule
$\phi_*(\alpha)(\psi):=\alpha\wedge\phi^*(\psi)(1)$. Le membre de droite est bien d\'efini, car le courant
$\alpha\wedge\phi^*(\psi)$ est \`a support compact, et peut donc s'\'etendre de mani\`ere unique en forme lin\'eaire
continue sur les formes lisses (avec la topologie d\'efinie par les semi-normes $p_{m,K}$).

 Si 
$r<p$ ou $s<p$, l'image directe est nulle. De plus, comme $\phi$ est analytique, $\phi_*$ commute avec
l'op\'erateur $\overline\partial$ sur les courants \`a support $\phi-$propre.
En effet
$\phi_*(\overline\partial\alpha)(\psi)=\overline\partial\alpha(\phi^*(\psi))=\alpha(\overline\partial\phi^*(\psi))$
soit comme $\overline\partial\phi^*(\psi)=\phi^*(\overline\partial \psi)$, $\overline\partial(\phi_*(\alpha))(\psi)$.

  Supposons maintenant $Y:=X\times F$, avec $F$ une vari\'et\'e analytique de dimension $p$, $X$ une vari\'et\'e analytique
de dimension $n$, et
$\phi:Y\to X$ la projection standard.

On d\'efinit les formes {\it verticales} sur $Y$, de la mani\`ere suivante.
Une forme $\psi$ est verticale si, \'etant donn\'e des coordonn\'ees locales
$y_1,\dots,y_p$ sur un ouvert $V$
$F$ et des coordonn\'ees locales $x_1,\dots,x_n$ sur un ouvert $U$
$X$, la forme s'\'ecrit localement sur $U\times V$ sous la forme $\sum_{I,J}{\psi_{I,J} dy^I\wedge d{\overline y}^J}$, avec
$\psi_{I,J}$ des fonctions sur $U\times V$; on pose comme d'habitude $dy^I:=dy_{i_1}\wedge\dots\wedge 
dy_{i_q}$. En particulier, les indices de degr\'e d'une telle forme sont inf\'erieurs \`a
$p$.

 Soit maintenant $\psi$ une forme sur $Y$ de bidegr\'e $(r,s)$ arbitraire. 
On se place au-dessus d'un ouvert $U$ de $X$, avec des coordonn\'ees $x_1,\dots,x_n$.
 $\psi$ s'\'ecrit de mani\`ere unique sous la forme :
$\psi=\sum_{I,J\subset\{1,\dots,n\}}{\psi_{I,J}{dx^I\wedge 
d\overline{x}^J}}$, o\`u les $\psi_{I,J}$ sont des formes verticales. Comme les indices de degr\'e de $\psi_{I,J}$
sont inf\'erieurs \`a $p$, les multiindices $I$ et $J$ intervenant dans la somme sont n\'ecessairement de
longueur respectivement sup\'erieures \`a $r':=r-p$ et $s':=s-p$.
 Si $r<p$, ou $s<p$, on a vu que $\phi_*([\psi])=0$. Supposons $r,s\ge p$.
 
$$\psi=\sum_{I=(i_1,\dots,i_{r'}),i_1<\dots<i_{r'},J=(j_1,\dots,j_{s'}),j_1<\dots<j_{s'}}{\psi_{I,J}dx^I\wedge 
d\overline{x}^J}+\psi',$$ 
o\`u $\psi'$ regroupe les termes correspondant \`a $card(I)>r'$ ou $card(J)>s'$. On a
vu  que 
$\phi_*([\psi'])=0$. 
On en d\'eduit :  $$\phi_*(\psi)=\sum_{I,J,card(I)=r',card(J)=s'}{\phi_*(\psi_{I,J}){dx^I\wedge
d\overline{x}^J}},$$  o\`u $\psi_{I,J}$ sont des formes verticales de bidegr\'e $(p,p)$.

Il reste donc \`a calculer $\phi_*(\psi)$, pour $\psi$ verticale de bidegr\'e $(p,p)$. Soit dans ce cas
$\psi_x$ la forme sur $F$ obtenue pour une valeur $x$ fix\'ee. Alors, d'apr\`es le th\'eor\`eme de Fubini,
$\mu(x):=\int_F{\psi_x}$ est {\it int\'egrable} et :
$$\int_{X\times F}{\psi\wedge f(x)dx\wedge d{\overline x}}=\int_X{f(x)dx\wedge d{\overline x}\int_F{\psi_x}}.$$
En particulier, si $\psi$ est int\'egrable, on voit donc que le courant $\phi_*(\psi)$ est repr\'esent\'e par la
fonction int\'egrable $\mu(x):=\int_F{\psi_x}$. Si $\psi$ est de classe ${\cal C}^k$ \`a support compact, on a (cf.
par exemple Schwartz,\cite{Schwartz}) $\mu(x):=\int_F{\psi_x}$ est de classe ${\cal C}^k$. On en d\'eduit, pour tous
les bidegr\'es :

    \begin{lem} Si $\psi$ est de classe ${\cal C}^k$, 
$\phi_*(\psi)$ est de classe 
${\cal C}^k$. 
 \end{lem} 

 Supposons $(r,s)=(N,p)$. Alors $\psi$ s'\'ecrit sur $U\times F$ de mani\`ere unique $\psi=\psi'\wedge dx$,
avec $dx:=dx_1\wedge\dots\wedge dx_n$ ($n:=\dim(X)$) et $\psi'$ une
$(p,p)-$forme verticale. On a vu que : $\phi_*(\psi)=\int_{F}{\psi'_x}dx$. Remarquons que $\psi'_x$ peut aussi
s'\'ecrire le r\'esidu de Poincar\'e-Leray 
$res_{F_x}{\psi/((X_1-{x_1})\dots(X_n-x_n))}$ ($F_x:=\phi^{-1}(x)$), 
avec pour $P\in Y$, $X_i(P):=x_i(\phi(P))$. 

Soit maintenant $\alpha$ un courant l. r. de bidegr\'e $(N,p)$ sur $Y$. 
Pour tout point $z\in Y$, on a donc un voisinage ouvert $U_{z}$ de $z$, tel que sur 
$U_z$ 
on ait : 
$\alpha=[\omega/g_z]\wedge{\overline\partial[1/{f_1}_z]}\wedge\dots\wedge{\overline\partial[1/{f_p}_z]}$, 
o\`u $\omega$ est une $N-$forme holomorphe, et $({f_1}_z,\dots,{f_p}_z,g_z)$ une suite 
r\'eguli\`ere 
de fonctions holomorphes sur $U_z$. 

On suppose que la restriction de $\phi$ au 
support $Z$
de 
$\alpha$ est propre. Donc la fibre $F_x=\phi^{-1}(x)$ coupe $Z$ en un 
nombre 
fini de points $P_i(x):=(x,y_i(x))$.

  De ce qui pr\'ec\`ede on d\'eduit le lemme suivant, qui nous servira par la 
suite: 

   \begin{lem} 
$\phi_*(\alpha)$ est une $n$-forme m\'eromorphe sur $X$; sur $U$, en coordonn\'ees locales
$x_1,\dots,x_n$, $\phi_*(\alpha)=\mu(x)dx$. On a, en dehors de l'hypersurface $\phi(Pol(\alpha))$ de $X$ :
 $$\mu(x)=Res(\alpha/ (X_1-x_1)\dots (X_n-x_n))(1),$$
 o\`u 
${Res(\alpha/h_1\dots h_n}$ est d\'efini localement par la formule:
$$Res(\alpha/(X_1-x_1)\dots (X_n-x_n):=Res{\omega/(f_1\dots f_p (X_1-x_1)\dots (X_n-x_n))}.$$
On peut encore \'ecrire :
$\mu(x)=\sum_{i}{Res_{P_i(x)}{\omega'_x/{f_{1}\dots
f_{p}} }}$, o\`u $\omega'$ est la 
$p-$forme verticale telle que $\omega=\omega'\wedge dx$, les $P_i(x)$ sont les points d'intersection de la fibre $F_x$ avec
le support de $\alpha$, et
$Res_{P_i(x)}$ est le r\'esidu ponctuel de Grothendieck.
  \end{lem} 

\dem Pla\c{c}ons-nous au-dessus d'un ouvert de $X':=X-\phi(Pol(\alpha))$, ayant des coordonn\'ees locales $x_1,\dots,x_n$.
Alors, \'ecrivons $\alpha:=\alpha'\wedge dx$, avec $\alpha'$ un {\it courant vertical} (un courant \'etant une forme \`a
coefficients distributions, on a une d\'efinition analogue \`a celle donn\'ee pour les formes.
Alors, $\alpha'$ est un courant vertical $\overline\partial-$ferm\'e de bidegr\'e $(p,p)$ \`a support $\phi-$propre. Donc :
$\phi_*(\alpha')$ est la fonction holomorphe d\'efinie par $\int_{F_x}{\alpha'}$.
Mais, d'apr\`es ce qu'on a vu, si $\alpha'$ est donn\'e localement aux points par $Res_{f_1,\dots,f_p}{\omega'/{f_1\dots
f_p}}$, alors :
$\int_{F_x}{\alpha}$ est la somme des r\'esidus ponctuels.
Donc $\phi_*(\alpha)$ est une $n-$forme holomorphe sur $X'$, donn\'ee localement par $\mu(x)dx$, et
$\mu(x)=\sum_{i}{Res_{P_i(x)}{ \omega'/f_1\dots f_p}}$. De plus, comme $\phi_*(\alpha)$ reste d'extension standard, c'est
un courant valeur principale, \`a p\^ole dans $\phi(Pol(\alpha))$.
  \CQFD 

On a calcul\'e l'image directe d'une forme lisse ou d'un courant l. r. \`a support $\phi-$propre pour une projection $\phi$. 

 Supposons maintenant qu'on ait comme morphisme analytique $\phi:Y\to X$ une {\it 
submersion}, d'un espace analytique r\'eduit $Y$ sur une vari\'et\'e analytique $X$ de dimension $n$. Cela signifie que
$\phi$ est localement trivialisable, i.e. qu'il existe pour tout point
$z\in Y$ un voisinage ouvert
$U_z$  et un biholomorphisme 
$\mu_z: D_{p(z)}\times V_z \simeq U_z$, tel que $\phi=\pi\circ\mu_z^{-1}$, avec $\pi$ 
la projection naturelle $\pi : D_{\phi(z)}\times V_z\to D_{\phi(z)}$, et $D_{\phi(z)}$ 
un voisinage ouvert de 
$\phi(z)$. Si $y\in Y$ est un point r\'egulier, cela revient \`a dire la diff\'erentielle $d\phi_y$ est de rang maximum
$n$, d'apr\`s le th\'eor\`eme des fonctions implicites.

Soit $\alpha$ un courant sur $Y$ tel que la 
restriction de $\phi$ sur le support de $\alpha$ est propre. 

\begin{lem} Si $\alpha=[\psi]$, avec $\psi$ une forme de classe ${\cal C}^k$, $\phi_*(\alpha)=[\rho]$, 
avec 
$\rho$ une forme de classe ${\cal C}^k$. 
\end{lem} 

   \dem 
On consid\`ere un recouvrement d'ouverts relativement compacts $(U_{i})_{i\in I}$ de $Y$ par 
des ouverts sur lesquels 
$\phi$ est trivialisable. On se donne une partition de l'unit\'e associ\'ee $\chi_i$. 

Soit 
$\psi$ une forme ${\cal C}^k$ sur $Y$, $\phi-$propre. On a : 
$\psi=\sum_{i\in I}{\rho\chi_i}$, et donc 
$\phi_*(\psi)=\sum_{i\in I}{\phi_*(\chi_i\psi)}$. 
La somme est bien d\'efinie, car pour tout point $x\in X$, il existe un voisinage 
$U_x$ relativement compact, tel que $\phi^{-1}(\overline{U_x})\cap supp(\psi)$ est 
compact, donc 
au-dessus de 
$U_x$ $\sum_{i\in I}{\psi\chi_i}$ est \'egale \`a une somme finie $\sum_{i\in J}{\psi\chi_i}$, et sur $U_x$ l'image
directe est la somme des images directes. Chaque terme de la somme est ${\cal C}^k$,  d'apr\`es ce qu'on a vu pour
les projections; donc la somme est de classe ${\cal C}^k$. 

\CQFD 

 On se donne maintenant un morphisme analytique $\phi :Y\to X$ entre deux espaces analytiques de 
dimension pure $X$ et $Y$, avec $\dim(Y)=\dim(X)+p$. On suppose qu'il existe un sous-ensemble analytique $S$
d'int\'erieur vide, telle que $\phi^{-1}(S)$ est aussi d'int\'erieur vide et $\phi:Y-\phi^{-1}(S)\to X-S$ est une
submersion.

 On consid\`ere un courant l. r. $\alpha$ sur $Y$, tel que la 
restriction de 
$\phi$ au support de 
$\alpha$ est propre. 

\begin{lem} Supposons $\alpha$ l. r. de bidegr\'e $(q+p,p)$ sur $Y$, d'ensemble 
polaire 
$Pol(\alpha)$. Alors $\phi_*(\alpha)$ est de la forme $[\omega]$, sur $X$, avec $\omega$ une $q-$forme m\'eromorphe
\`a  p\^ole 
contenu dans $\phi(Pol(\alpha))\cup Sing(X)$. Si $\alpha$ est 
$\overline\partial-$ferm\'e, 
$\omega$ est ab\'elienne sur $X$. 
\end{lem}

\dem
Premi\`erement, en dehors de $\phi(Pol(\alpha))\cup Sing(Y)$, $\phi_*(\alpha)$ est l'image 
directe d'un courant de bidegr\'e $(q+p,p)$ $\overline\partial-$ferm\'e. C'est donc un 
courant $\overline\partial-$ferm\'e de bidegr\'e $(q,0)$, donc une $q-$forme holomorphe. 

D'apr\`es le lemme \cite{HenPas}, cette forme holomorphe s'\'etend sur $T$ en une forme 
m\'eromorphe $\omega$. De plus, si $\alpha$ est $\overline\partial-$ferm\'e, alors 
$\phi_*(\alpha)$ est $\overline\partial-$ferm\'e; il en est donc de m\^eme de $[\omega]$. 

\CQFD 

\begin{rem} Supposons $\dim(X)=\dim(Y)$, et $\phi: Y\to X$ propre. 
Si $\alpha=[\psi]$ pour une forme m\'eromorphe 
$\psi$ sur $X$, la forme m\'eromorphe $\omega$ telle que $\phi_*([\psi])=[\omega]$ s'appelle 
la {\it trace} de 
$\omega$ pour $\phi$, et se note encore $\phi_*(\psi)$. 
Si $\psi$ est ab\'elienne, $\phi_*(\psi)$ l'est aussi.\end{rem}

\begin{lem} Supposons $S=\emptyset$, i.e. $\phi:Y\to X$ est une submersion. Alors, si $\alpha$ est l. r. de bidegr\'e
$(r,s)$,
\`a support
$\phi-$propre,
$\phi_*(\alpha)$ est l. r. de bidegr\'e $(r-p,s-p)$.
\end{lem}
\dem Soit $\beta:=\phi_*(\alpha)$.
 D'apr\`es le th\'eor\`eme de Remmert, si $Z$ est le support de $\alpha$, $Z':=\phi(Z)$, qui est le support de
$\beta$, est un sous-ensemble analytique de
$X$. Alors, $\overline{I_Z}$ annule
$\alpha$, donc $\overline{I_{Z'}}$ annule $\beta$.
  On v\'erifie \'egalement que $\beta$ est $\overline\partial-$ferm\'e sur un ouvert de Zariski dense de $Z'$.
Enfin, l'image directe d'un courant d'{\it extension standard} par une submersion reste d'extension standard.  Ces trois
propri\'et\'es suffisent d'apr\`es ce qu'on a vu pour caract\'eriser les
courants l. r. . On en d\'eduit que
$\phi_*(\alpha)$ est l. r. .
\CQFD

\subsection{Image inverse d'un courant par une submersion.}

Soit $\phi:X\to Y$ une submersion analytique; l'image directe d'une forme lisse \`a support compact est une forme
lisse
\`a support compact. On d\'efinit alors l'image inverse d'un courant $\beta$ sur $X$ par 
$(\phi^*(\beta),\psi):=(\beta,\phi_*(\psi))$. On v\'erifie $\phi^*(\beta\wedge\psi)=\phi^*(\beta)\wedge\phi^*(\psi)$ pour
une forme lisse $\psi$. De plus :

 \begin{pro} Si $\beta$ est un courant l. r., 
 $\phi^*(\beta)$ l'est aussi, et le r\'esidu commute avec l'image inverse : 
$$\phi^*([1/g]\wedge{\overline\partial 
[1/f_1]}\wedge\dots\wedge{\overline\partial[1/f_k]})=[1/\phi^*(g)]{\overline\partial[1/\phi^*(f_1)]}
\wedge\dots\wedge{\overline\partial[1/\phi^*(f_k)]}.$$ 
  \end{pro}

 \dem D'abord, $\phi^*(f_1),\dots,\phi^*(f_k),\phi_*(g)$ 
forment toujours une suite r\'eguli\`ere, puisque la codimension de 
$\phi^*(f_1)=\dots,\phi^*(f_k)=\phi_*(g)=0$ est la m\^eme que celle de $f_1=\dots=f_k=g=0$. 
Il suffit de v\'erifier la proposition lorsque $\phi$ est une projection, puisqu'une submersion est localement
\'equivalente \`a une projection.
Soit $\psi$ une forme lisse \`a support compact sur $Y=F\times X$.
On obtient :
$$
\int_{F\times X}{[1\phi^*(g)]{\overline\partial[1/\phi^*(f_1)]}\wedge\dots\wedge
{\overline\partial[1/\phi^*(f_p)]}\wedge\psi}=\int_X{ [1/g]\wedge\overline\partial
{[1/f_1]}\wedge\dots\wedge\overline\partial[1/f_p]\wedge\int_{F_x}{\psi}},$$ 
en \'ecrivant les d\'efinitions, par la formule de Fubini g\'en\'eralis\'ee.
\CQFD 

On a de plus :

\begin{lem} Soit $\phi:Y\to X$ une submersion. Soit $Z$ un sous-ensemble analytique de $X$.
Soit $\omega$ une forme m\'eromorphe sur $Z$. Alors : $\phi^*(\omega\wedge[Z])=\phi^*(\omega)\wedge[\phi^{-1}(Z)]$.
\end{lem}

\subsection{Transformation d'Abel-Radon.} 

On se donne les objets suivants : 

i) Un espace analytique r\'eduit $X$ de dimension pure $N$. 

ii) Un espace de param\`etres $T$, espace analytique r\'eduit de dimension pure. 

iii) Une vari\'et\'e d'incidence $I\subset T\times X$, sous-ensemble analytique, avec les projections naturelles $p_1: I\to
T$ et
$p_2:I\to X$  et les propri\'et\'es suivantes: 

1) $I$ est de dimension pure $\dim(T)+p$.

2) La projection $p_1:I\to T$ est propre, et telle qu'il existe un sous-ensemble analytique $S\subset T$
d'int\'erieur vide, tel que ${p_1}^{-1}(S)$ est encore d'int\'erieur vide, et que la restriction de $p_1$ \`a
$I-{p_1}^{-1}(S)$ est une submersion. 

3) La projection $p_2:I\to X$ est une submersion. 

Cette troisi\`eme condition nous permet de d\'efinir, comme on l'a vu ci-dessus, 
l'image inverse ${p_2}^{-1}(\alpha)$ d'un courant $\alpha$ sur $X$. Comme $p_1$ est propre, 
on peut prendre l'image directe de ce courant par $p_1$, qui est par d\'efinition la {\it 
transform\'ee d'Abel-Radon} de $\alpha$: 

 \begin{defi}$${\cal R}(\alpha)={p_1}_*({p_2}^*(\alpha)).$$
\end{defi} 
 
 La transformation associe donc \`a un courant de bidegr\'e $(r,s)$ un courant de bidegr\'e $(r-p,s-p)$. 

Soit $Y$ un sous-ensemble analytique de $X$ de codimension pure $p$. Alors, posons $I_Y:=p_2^{-1}(Y)\subset I$.
$I_Y$ d\'efinit une famille de z\'eros-cycles sur $Y$, avec les deux projections $\pi_1:I_Y\to T$ et $\pi_2:I_Y\to
Y$. On peut d\'efinir la transformation d'Abel d'une forme m\'eromorphe $\omega$ sur $Y$, par ${\cal
A}(\omega):={\pi_1}_*(\pi_2^*(\omega))$, qui est une forme m\'eromorphe sur $T$ (cf.\cite{Fabre2}). 
Alors :
\begin{pro} ${\cal A}(\omega)={\cal R}(\omega\wedge[Y]).$
\end{pro}
\dem
On a ${p_2}^*(\omega\wedge[Y])=p_2^*(\omega)\wedge[I_Y])$.
De plus, ${p_1}_*(\psi\wedge[I_Y])={{p_1}_{\vline I_Y}}_*(\psi)$, avec $\psi:=p_2^*(\omega)$. 
\CQFD
 \begin{lem} $\cal R$ commute avec les op\'erateurs $d,\overline\partial,\partial$.
\end{lem} 

 \dem 
 On a d\'ej\`a vu que les op\'erateurs $d,\overline\partial,\partial=d-\overline\partial$ sur les courants commutent
avec l'image directe par  une application holomorphe propre. 
On a vu aussi que ces op\'erateurs sur les courants commutent avec l'image inverse par une 
submersion analytique. 
\CQFD 

 \begin{cor} 
Supposons $\alpha$ une courant $\overline\partial-$ferm\'e de bidegr\'e $(q+p,p)$. Alors 
${\cal R}(\omega)$ d\'efinit une $q-$forme ab\'elienne sur $T$.\end{cor} 

\begin{pro} 
Supposons $p_1:I\to T$ une submersion. 
 Si $\alpha$ est un courant l. r., alors ${\cal R}(\alpha)$ est encore un 
courant l. r. .\end{pro} 

 \dem 
On a vu que l'image inverse par une submersion d'un courant l. r. reste localement 
r\'esiduel, de m\^eme que pour l'image directe.
\CQFD 

    \subsection{Un th\'eor\`eme d'extension.} 
  Soit $D$ un domaine d'une vari\'et\'e analytique $X$. On consid\`ere la projection standard $\phi:D\times \C \to D$.

 On se 
donne sur $U:=D\times \C$ un courant $\alpha$ {\it r\'esiduel} de bidegr\'e $(n+1,1)$, $\overline\partial-$ferm\'e.
 On a donc une $(n+1)-$forme m\'eromorphe $\psi$ sur $U$, telle que $\alpha=\overline\partial[\psi]$. On suppose que
la restriction de
$\phi$  au support de $\alpha$ est propre. 
  Comme la restriction de $\phi$ au support de $\alpha$ est propre, 
on peut d\'efinir le courant image directe $\phi_*(\alpha)$. 

Sur un ouvert $U$ de $X$ avec des coordonn\'ees locales $x_1,\dots,x_n$, on note
$dx:=dx_1\wedge\dots\wedge dx_n$, et $\overline{dx}:=\overline{dx_1}\wedge\dots\wedge 
\overline{dx_n}$. 

On d\'efinit la $1-$forme m\'eromorphe verticale $\omega$ telle que
$\psi=\omega\wedge dx$; on note $\omega_z$ sa restriction sur $F_z:=\phi^{-1}(z)$ Alors, ce qui pr\'ec\`ede montre que :

 \begin{lem} 
\label{expli} Sur $U$, $\phi_*(\alpha y^k)=u_k(x)dx$, avec $u_k$ une fonction holomorphe, d\'efinie
par : $u_k(x)=\sum{Res_{y_i(x)}{\omega_x}}$, 
 o\`u les points $y_i(x)$ sont les p\^oles de $\omega_x$ sur $F_x$, et o\`u $Res$ d\'esigne le r\'esidu classique d'une
$1-$forme 
m\'eromorphe dans $\C$. 
  \end{lem}

\begin{lem} Supposons que pour tout $k\ge 0$, on ait $u_k(x)=0$. Alors, 
$\alpha=0$. \label{zero}\end{lem} 

\dem 
En effet, supposons que $\alpha$ soit non nul. Alors $\psi$ est m\'eromorphe, et pour un $x$ g\'en\'erique,
$\omega_x$ est aussi m\'eromorphe. On ne peut donc pas avoir :
$\sum_{y_i(x)}{res_{y_i(x)}{\omega_x y^k}}=0$ 
 pour tout $k$. 
\CQFD 

Consid\'erons le corps $K$ des fonctions m\'eromorphes sur $D$. On peut multiplier le courant l.r. $\alpha$, \`a
support
$\phi-$propre, par n'importe quelle fonction m\'eromorphe $Q\in K[y]$, puisque le p\^ole de $Q$ coupe le support de
$\alpha$ proprement. L'ensemble des
$Q$ tels que
$Q\alpha=0$ forme un id\'eal dans l'anneau principal $K[y]$, est donc engendr\'e par un \'el\'ement
$P=y^d+a_1(x)y^{d-1}+\dots+a_d(x)$, de degr\'e minimal $d$.

\begin{lem} Les fonction m\'eromorphes $a_i(x)$ sont en fait holomorphes.
\end{lem}
\dem En effet, soit $D'\subset D$ l'ouvert de Zariski sur lequel toutes les $a_i$ sont holomorphes. Au-dessus de $D'$,
on a $Supp(\alpha)=\{P=0\}$. Supposons
$D'\not=D$. Alors, soit $x\in D-D'$, et $U_x$ un voisinage ouvert de $x$ relativement compact dans $D$. Alors
$\phi^{-1}(\overline{U_x})\cap Supp(\alpha)$ ne serait pas compact, ce qui contredit l'hypoth\`ese. 
\CQFD

En prenant l'image directe par $\phi$ de $y^kP(x,y)\alpha=0$, on obtient :
\begin{lem}
On a $u_{k+d}(x)+a_1(x)u_{k+d-1}(x)+\dots+a_d(x)u_{k}(x)=0$ pour tout entier $k\ge 0$.
\end{lem}

Soit la fonction holomorphe
$r(x,y):=y^{d-1}u_0(x)+y^{d-2}(u_1(x)+a_1(x)u_0(x))+\dots+(u_{d-1}(x)+a_1(x)u_{d-2}(x)+\dots+a_{d-1}(x)u_0(x))$.

Alors : 

\begin{lem} $\alpha=\overline\partial[r(x,y)dx\wedge dy/P(x,y)]$.\end{lem} 

\dem 
Consid\'erons la s\'erie formelle $\psi:=\sum_{k\ge 0}{u_k(x)/y^{k+1}}$. D'apr\`es le lemme pr\'ec\'edent, on a
l'\'egalit\'e formelle :
$P(x,y)\psi=r(x,y)$. D'autre part, consid\'erons la fonction m\'eromorphe : $r(x,y)/P(x,y)$. Pour $x$ fix\'e, on
peut r\'eduire la fraction rationnelle en \'el\'ements simples, et d\'evelopper en puissances de $1/y$. On obtient
donc une s\'erie, qui converge uniform\'ement pour $\vline y\vline\ge R(x)> \max_i{y_i(x)}$. Les coefficients
obtenus doivent
\^etre
\'egaux aux $u_k(x)$.

 Soit $\beta:=\overline\partial[r(x,y)dx\wedge dy/P(x,y)]$. 
D'apr\`es le lemme \ref{expli} pr\'ec\'edent, on a $\phi_*(\beta y^k)=v_k(x)dx$, avec 
$v_k(x):=\sum_{y_i(x)}{res_{y_i(x)}{r(x,y)y^k/P(x,y)dy}}$ pour tout $k\ge 0$, soit encore, pour un $R(x)$ assez
grand :
$$1/2i\pi \int_{\vline y\vline = R(x)}{r(x,y)y^k/P(x,y)dy}.$$ L'int\'egrande vaut encore, comme on vient de voir,
$\sum_{j\ge  0}{u_k(x)/y^{j+1-k}}$, et la somme converge uniform\'ement pour $\vline y\vline = R(x)$. L'int\'egrale
de la somme est donc la somme des  int\'egrales; 
un seul terme est non nul, c'est $u_k(x)$. On a donc finalement $\phi_*(\beta 
y^k)=u_k(x)dx$, pour tout $k\ge 0$. 
D'apr\`es le lemme \ref{zero}, on obtient $\alpha=\beta$. 
\CQFD 

\begin{lem} Le d\'eterminant de la matrice $M=(u_{d+i-j-1}(x))_{1\le i,j\le 
d}$ est non identiquement nul.\end{lem} 

\dem Supposons que le d\'eterminant 
soit identiquement nul. On consid\`ere la matrice comme un endomorphisme 
de $K^d$, o\`u $K$ est le corps des fonctions m\'eromorphes sur $D$. 
Si le d\'eterminant est nul, on a un vecteur propre associ\'e \`a la valeur 
propre nulle. On a donc une relation lin\'eaire 
$b_1(x)u_{d-1+k}(x)+\dots+b_d(x)u_k(x)=0$, pour $k=0,\dots,d-1$. Mais 
alors, les \'equations 
$u_{k+d}(x)+a_1(x)u_{k+d-1}(x)+\dots+a_d(x)u_{k}(x)=0,$ valables pour tout 
$k$, impliquent que $b_1(x)u_{d-1+k}(x)+\dots+b_d(x)u_k(x)=0$ reste valide 
pour $k\ge d$. Mais cela implique 
$$\phi_*(y^k(b_1(x)y^{d-1}+\dots+b_d(x))\alpha)=0,$$ pour tout $k\ge 0$. On en 
d\'eduit que $(b_1(x)y^{d-1}+\dots+b_d(x))\alpha=0$.
Cela implique que $(b_1(x)y^{d-1}(x)+\dots+b_d(x))$ est multiple de $P(x,y)$, ce 
qui est impossible, car le polyn\^ome en $y$ $P(x,y)$ est de degr\'e 
$d$.\CQFD 

\begin{rem} Dans le cas o\`u le support de $\alpha$ est "r\'eduit", 
 on peut calculer le d\'eterminant de la mani\`ere suivante. Supposons
$\alpha:=\omega\wedge[Y]$, o\`u $Y$ est le support de 
$\alpha$ et $\omega$ une $n-$forme m\'eromorphe sur $Y$. On 
peut \'ecrire $\omega_i=f_i(x)dx$, sur les diff\'erentes composantes $Y_i$ de $Y$ au-dessus 
d'un voisinage ouvert de $x_0\in D$. Alors, le d\'eterminant s'\'ecrit 
$$\Pi_{i<j}{(y_i(x)-y_j(x))^2}\Pi_{i=1}^d{f_i(x)}.$$ 
\end{rem}

\begin{thm} 
 Supposons que dans la situation pr\'ec\'edente, tous les $u_k(x)$ se 
prolongent m\'eromorphiquement \`a un domaine $\tilde D$ contenant $D$. 
Alors, $\alpha$ se prolonge \`a un courant l. r. $\tilde\alpha$ d\'efini 
sur ${\tilde D}\times \C$, et sans composante verticale, i.e. son support, ne contient pas
d'hypersurface de la forme $\phi^{-1}(H),H\subset {\tilde D}$. 
\end{thm} 

\dem On consid\`ere les \'equations 
$u_{k+d}(x)+{\tilde a_1}u_{k+d-1}(x)+\dots+{\tilde a_d}u_{k}(x)=0$, soit
encore:
${\tilde a_1}u_{k+d-1}(x)+\dots+{\tilde a_d}u_{k}(x)=-u_{k+d}(x)$, o\`u on consid\`ere les ${\tilde a_i}$ comme des
inconnues. On \'ecrit les \'equations pour
$k=0,\dots,d-1$. On obtient un syst\`eme carr\'e d'\'equations
lin\'eaires, \`a coefficients dans le corps $K$ des fonctions m\'eromorphes sur $\tilde D$.

Le lemme 
ci-dessus montre que le d\'eterminant du syst\`eme est non identiquement nul sur $D$, donc, comme 
il est m\'eromorphe sur $\tilde D$, il n'est pas non plus identiquement nul sur $\tilde D$. Il admet donc des
solutions
$\tilde{a_i}$ uniques, qui doivent donc co\"\i ncider avec $a_i$ sur $D$ (qui sont les solutions uniques du
syst\`eme sur
$D$). 

On prolonge donc $P(x,y)$, ainsi que $r(x,y)$, en fonctions m\'eromorphes ${\tilde P}$ et $\tilde r$ sur ${\tilde
D}\times
\C$. Consid\'erons le courant $\overline\partial[{\tilde r}(x,y)/{\tilde P}(x,y)dx\wedge dy]$. 
C'est un courant r\'esiduel sur ${\tilde D}\times \C$, qui co\"\i ncide avec $\alpha$ sur $D\times \C$, 
et qui est sans composante verticale.
\CQFD 

Remarquons la possibilit\'e, m\^eme si les $u_k(x)$ se prolongent en fonctions {\it holomorphes}, que les prolongements
${\tilde a}_i(x)$ sur ${\tilde D}$ des coefficients $a_i(x)$ holomorphes sur $D$ soient {\it m\'eromorphes} sur
${\tilde D}$.

Le lemme suivant est une variante du lemme \ref{HenPas} ci-dessus, variante qui nous sera utile par la suite.

\begin{lem}\label{mero} Soit $Y$ une hypersurface analytique de $D\times\C$, telle que la restriction de
$\phi:D\times\C\to D$
\`a $Y$ soit propre. Soit $Z$ un sous-ensemble analytique de $Y$ d'int\'erieur vide dans $Y$.
$\alpha$ un courant l.r. de bidegr\'e
$(n+1,1)$ sur
$D\times\C-Z$, \`a support dans $Y-Z$, $\overline\partial-$ferm\'e. Si pour tout $k\ge 0$, $u_k(x)$ se prolonge
m\'eromorphiquement
\`a travers
$Z':=\phi(Z)$, alors
$\alpha$ se prolonge de mani\`ere unique en un courant l. r. $\beta$ \`a travers $Z$. Si les $u_k$ se prolongent
holomorphiquement, $\beta$ est $\overline\partial-$ferm\'e. 
\end{lem}

\dem On associe \`a $\alpha$ comme ci-dessus le polyn\^ome $P(x,y):=y^d+a_1(x)y^{d-1}+\dots+a_d(x)$, dont les
coefficients sont holomorphes sur $D-Z'$. Le support de $\alpha$ dans $D-Z$, d'\'equation $\{P(x,y)=0\}$, est
contenu dans
$Y$. Donc son adh\'erence dans $D\times \C$ est une r\'eunion de composantes de $Y$. On en d\'eduit que les
$a_i(x)$ se prolongent {\it holomorphiquement} dans $D$. On d\'efinit \'egalement comme ci-dessus la
fonction
$r(x,y):=y^{d-1}u_0(x)+y^{d-2}(u_1(x)+a_1(x)u_0(x))+\dots+(u_{d-1}(x)+a_1(x)u_{d-2}(x)+\dots+a_{d-1}(x)u_0(x))$,
avec $\alpha=\overline\partial[r(x,y)/P(x,y)dx\wedge dy]$ sur $(D-Z')\times\C$.  Si les $u_0,\dots,u_{d-1}$ se
prolongent holomorphiquement, $r(x,y)$ d'apr\`es son expression se prolonge aussi holomorphiquement, et
$\tilde\alpha:=\overline\partial[{\tilde r} dx\wedge dy/{\tilde P}]$ est un prolongement
$\overline\partial-$ferm\'e de $\alpha$ sur $D\times\C$, \`a support dans $Y$, et co\"\i ncide avec $\alpha$ en
dehors de $Z$. Supposons que les $u_k$ se prolongent {\it m\'eromorphiquement} \`a travers $Z'$.
Soit $x_0\in Z'$.
Soit $G(x)$ un fonction holomorphe en $x_0$ telle que $Gu_k(0\le k\le d-1)$ soient holomorphes dans un
voisinage ouvert $U_{x_0}$ de $x_0$. 
Alors $Gr$ se prolonge holomorphiquement, soit $R(x,y)$ sur $U_{x_0}\times\C$. Alors
$1/G\overline\partial[R/{\tilde P}]$ est un prolongement l.r. de $\alpha$ sur $U_{x_0}\times\C$, dont le support est
contenu dans $Y$, et donc co\"\i ncide avec $\alpha$ en dehors de $Z$. Le prolongement est unique, puisque $Z$ est de
codimension deux.
\CQFD

 \subsection{Un cas particulier : la transformation par rapport aux $p-$plans.} 

Soit $G(1,N)$ la grassmanienne des $p-$plans complexes, et $I_{\P_N}\subset G(p,N)\times \P_N$ la vari\'et\'e
d'incidence; on note $p_1$ et $p_2$ les deux projections standard $p_1:I_{\P_N}\to G(p,N)$ et
$p_2:I_{\P_N}\to \P_N$. Remarquons que $I_{\P_N}$ est lisse.

 Soit $U$ un ouvert lin\'eairement $p-$concave de $\P_N$, i.e. r\'eunion de $p-$plans
complexes (pour $p=N-1$, on dit simplement : lin\'eairement concave). On d\'efinit l'ouvert $U^*\subset G(p,N)$
comme l'ouvert correspondant aux
$p-$plans contenus dans
$U$, et la vari\'et\'e d'incidence $I_U=p_1^{-1}(U^*)\subset U^*\times U$; on note encore $p_1$ et $p_2$ les deux
projections standard
$p_1:I_U\to U^*$ et
$p_2:I_U\to U$, ou $p_1^U$ et $p_2^U$ lorsqu'on a besoin de pr\'eciser l'ouvert.

On v\'erifie : 

i) $\dim(I_U)=\dim(U^*)+p$, la fibre $p_1^{-1}(t)=\{t\}\times H_t$ \'etant de dimension $p$;

ii) $p_1:I_U\to U^*$ est propre, puisque $p_1^{-1}(K)={{p_1}^{\small\P_N}}^{-1}(K)$;

iii) $p_2:I_U\to U$ est une submersion. En effet, soit $(H_P,P)\in I_U$ un point quelconque, avec $P\in U$ et
$H_P\subset U$ un
$p-$plan passant par $P$. On se donne dans $H_P$ un $(p-1)-$plan $H'$ en dehors de $P$. On d\'efinit, sur un
voisinage ouvert
$U_P$ de
$P$ dans $U$, un morphisme analytique : $s_P:U_P\to I_U$, avec $p_2\circ s_P=Id_{U_P}$, de la mani\`ere suivante :
pour $Q\in U_P$, $s_P(Q)$ est le $p-$plan contenant $H'$ et passant par $Q$. Alors, pour tout vecteur tangent $v$
\`a $U$ en $P$, $v=dp_2(ds_P(v))$, avec $ds_P(v)$ un vecteur tangent \`a $I_U$ en $(H_P,P)$;
$dp_2:T_{(H_P,P),I_U}\to T_{P,U}$ est donc surjective, ce qui \'equivaut, par le th\'er\`eme des fonctions
implicites, comme $(H_P,P)$ est un point r\'egulier, \`a dire que $p_2$ est localement trivialisable en $(H_P,P)$.
On a donc une submersion analytique $p_2:I_U\to U$. 

Par l'axiome du choix, il existe une section $s$ de $p_2$, i.e. une application $s : U\to I_U$ telle que $p_2\circ
s=Id_U$. Mais il n'existe pas en g\'en\'eral de section {\it continue}. 

\begin{defi} On dit
que
$U$ est {\it contin\^ument} lin\'eairement $p-$concave, si on peut trouver une section $s:U\to I_U$ continue.
\end{defi} 

 \begin{lem} \label{GinHen}
Supposons que $U$ est contin\^ument lin\'eairement $p-$concave. Alors,
 soit $\alpha$ un courant $\overline\partial-$ferm\'e de bidegr\'e $(N,p)$.
$\alpha$ est $\overline\partial-$exact ssi ${\cal R}(\alpha)=0$.
\end{lem}

 \dem 
 La preuve pour les formes lisses d\'ecoule des formules de repr\'esentation int\'egrale 
de Henkin-Gindikin (cf. \cite{GHenkin}). Soit $\alpha$ un courant de bidegr\'e $(N,p)$. Alors,
$\alpha=\overline\partial\beta+\psi$, avec $\beta$ un courant et $\phi$ une forme lisse de bidegr\'e $(N,p)$. Alors
${\cal R}(\alpha)=0$ implique \`a ${\cal R}(\psi)=0$ donc $\phi=\overline\partial\mu$, donc
$\alpha=\overline\partial(\beta+\mu)$.
\CQFD

\begin{lem} Soit $f$ une fonction m\'eromorphe sur $\Delta^N$, o\`u $\Delta$ est le disque $\{z\in \C,\vline z\vline
<1\}$. Alors pour $(z_1^0,\dots,z_{i-1}^0,z_{i+1}^0,\dots,z_N^0)$ fix\'es dans un ouvert Zariski-dense de
$\Delta^{N-1}$, on peut consid\'erer $f_i(z_i):=(z_1^0,\dots,z_{i-1}^0,z_i,z_{i+1}^0,\dots,z_N^0)$ m\'eromorphe par
rapport \`a la variable restante $z_i\in \Delta$.
On suppose $f_i$ rationnelle lorsqu'elle est d\'efinie, pour chaque $i,1\le i\le N$. Alors 
$f$ est rationnelle. 
\end{lem} 

\dem 
 Le lemme a \'et\'e montr\'e pour $N=2$ par W. Kneser (\cite{Kneser}). La d\'emonstration s'adapte 
pour $N$ plus grand. 
\CQFD 

\begin{lem} \label{pseudoconcave} Soit $U$ un domaine de $\P_N$. Si $U$ contient 
une droite complexe, alors toute 
$q-$forme m\'eromorphe $\psi$ d\'efinie sur 
$U$ s'\'etend en une $q-$forme rationnelle sur $\P_N$. 
\end{lem} 

\dem 
On se ram\`ene ais\'ement au cas d'une fonction m\'eromorphe, en \'ecrivant $\psi$ 
comme combinaison, \`a coefficients m\'eromorphes, de formes rationnelles 
\'el\'ementaires $dx_{i_1}\wedge\dots\wedge dx_{i_q}$, avec $x_1,\dots,x_N$ un choix
convenable de coordonn\'ees affines. 

Le cas de la fonction m\'eromorphe peut se montrer de plusieurs mani\`eres. Une est 
donn\'ee dans \cite{Griffiths}. Une autre est bas\'ee sur le fait que l'ouvert $U$ est pseudoconcave au sens
d'Andreotti; on reverra cette notion plus loin. La preuve la plus \'el\'ementaire est bas\'ee sur le lemme
pr\'ec\'edent.  Soit un point
$P$ sur la droite complexe
$\Delta$ dont on suppose l'existence, et un hyperplan $H$ ne contenant pas
$P$. Choisir un syst\`eme de coordonn\'ees affines ayant $P$ comme origine et $H$ comme hyperplan \`a l'infini
revient \`a se donner $N$ points dans $H$; \`a savoir, les $N$ points $P_i (1\le i\le N)$, intersection de $H$ avec
les droites
$x_1=Cst,\dots,x_{i-1}=Cst,x_{i+1}=Cst,\dots,x_N=Cst$. Choisissons ces $N$ points dans $H\cap U$, et de sorte que
les $N$ droites $x_1=0,\dots,x_{i-1}=0,x_{i+1}=0,\dots,x_N=0$ soient contenues dans $U$ (il suffit de les choisir
suffisamment proches de $\Delta\cap H$). Alors,
$f$ est m\'eromorphe de $x_1,\dots,x_N$ dans un voisinage de l'origine, et elle est m\'eromorphe sur toutes les
droites
$x_1=x_1^0,\dots,x_{i-1}=x_{i-1}^0,x_{i+1}=x_{i+1}^0,\dots,x_N=x_N^0,$ pour tout
$(x_1^0,\dots,x_{i-1}^0,x_{i+1}^0,\dots,x_N^0)$ dans un ouvert de Zariski dense d'un voisinage ouvert de $0\in
\C^{N+1}$; donc rationnelle. $f$ s'exprime donc rationnellement en fonction des coordonn\'ees affines
$(x_1,\dots,x_N)$.
\CQFD

\begin{thm} Soit $U\subset \P_N$ un domaine lin\'eairement $1-$concave. Soit
$\alpha$ un courant de bidegr\'e
$(q+1,1)$, localement r\'esiduel, sur $U$. Supposons ${\cal R}(\alpha)=0$. Alors, $\alpha$ 
se prolonge de mani\`ere unique \`a un courant l. r. $\tilde\alpha$ sur $\P_N$;
$\tilde\alpha=\overline\partial[\psi]$, pour une forme rationnelle $\psi$ sur $\P_N$. 
\end{thm} 

\dem 
Tout d'abord, 
supposons que $U$ soit {\it contin\^ument} lin\'eairement $1-$concave, et $\alpha$ $\overline\partial-$ferm\'e. 
D'apr\`es le lemme \ref{GinHen}, on d\'eduit de ${\cal R}(\alpha)=0$ que $\alpha=\overline\partial \beta$, pour un
certain courant $\beta$ de bidegr\'e $(q+1,0)$. 
 Alors $\beta$ est d\'efini par une $(q+1)-$forme holomorphe en dehors du support $Y$ de 
$\alpha$. De plus, cette forme holomorphe doit se prolonger en une forme m\'eromorphe 
sur $U$, d'apr\`es le lemme \ref{HenPas}. On note cette forme m\'eromorphe $\psi$. 
On a alors: $\alpha=\overline\partial[\psi]+(\overline\partial{\beta-[\psi]})$, 
d'une part, et d'autre part sur chaque ouvert $U_i$ d'un recouvrement, 
$\alpha=\overline\partial[\psi_i]$, car $\alpha$ est localement r\'esiduel. 
Sur $U_i$, on a donc $\overline\partial[\psi_i-\psi]=\overline\partial{\beta-[\psi]}$; en particulier,
$[\psi_i-\psi+h]=(\beta-[\psi])$ pour une $(q+1)-$forme holomorphe $h$ sur $U_i$.
Mais le courant "valeur principale" $[\psi_i-\psi+h]$ est d'extension standard, le membre de droite nous montre
qu'il est nul en dehors de
$Y$; il est donc nul sur $U_i$.
On a donc $\alpha=\overline\partial[\psi_i]=\overline\partial[\psi]$ sur $U_i$. 

D'autre part, on a vu que la forme m\'eromorphe $\psi$ sur $U$ se prolonge en une forme rationnelle $\tilde\psi$ sur
$\P_N$. 
$\tilde\alpha:=\overline\partial[\tilde\psi]$ est donc un prolongement de $\alpha$. 

Dans le cas g\'en\'eral, consid\'erons une droite $\Delta$ contenue dans $U$, et un voisiange $U_\Delta$ de $\Delta$
contin\^ument lin\'eairement $1-$concave. Alors, on obtient un prolongement ${\tilde \alpha}_\Delta$ de $\alpha_{\vline
U_\Delta}$. Mais ${\tilde \alpha}_\Delta=\alpha$ sur $U$. En effet, si la diff\'erence \'etait non nulle, son support, qui
est un hypersurface analytique de $U$, devrait rencontrer la droite $\Delta$.

Pour la m\^eme raison, le prolongement est unique.
\CQFD

\begin{thm} 
Soit $U$ un ouvert de $\P_N$ lin\'eairement $p-$concave. 
Soit $\alpha$ un courant l. r. de bidegr\'e $(N,1)$ dans 
$U$. Alors ${\cal R}(\alpha)$ est m\'eromorphe sur $U^*$. $\alpha$ est $\overline\partial-$ferm\'e ssi
 ${\cal R}(\alpha)$ est
holomorphe.
 On suppose que ${\cal R}(\alpha)$ se prolonge m\'eromorphiquement dans un domaine 
${\tilde U}^*$. Alors $\alpha$ se prolonge dans $\tilde U$ comme courant l.
r. . 
\end{thm} 

Soit $U\subset {\tilde U}$ deux domaines lin\'eairement $1-$concaves. Soit $a\in U$. On note $\P_{N-1}^a\subset
G(1,N)$ l'ensemble des droites passant par $a$, et 
$D_a:= U^*\cap \P_{N-1}^a$, et ${\tilde D}_a:={\tilde 
U}^*\cap \P_{N-1,a}$. Soit $U_a:=\cup_{t\in D_a}{\Delta_t}-\{a\}$, ${\tilde U}_a:=\cup_{t\in
{\tilde D}_a}{\Delta_t}-\{a\}$. On d\'efinit, pour tout point $a\in {\tilde U}$, la projection
$p_a:{\tilde U}_a\to {\tilde D}_a$, qui \`a un point $P$ associe la droite qui la relie \`a $a$. 

 Soit $\alpha$ un courant l. r. \`a support dans l'ouvert lin\'eairement
$1-$concave
$U$, de support
$Y$. On se donne un point $a$ en dehors du support de $\alpha$. Alors $p_a:{U}_a\to {D}_a$, restreinte au support
de $\alpha$, est propre, et on peut d\'efinir son image directe ${p_a}_*(\alpha)$. 

\begin{lem} 
${p_a}_*(\alpha)={\cal R}(\alpha)_{\vline D_a}$
\end{lem}
\dem
On a ${\cal R}(\alpha)_{\vline D_a}={p_1}_*(p_2^*(\alpha)_{\vline I_a})$, avec $I_a:=p_1^{-1}(D_a)$. Mais
$p_2:I_a\to U_a$ est inversible, d'inverse $s_a:U_a\to I_a, x\mapsto (p_a(x),x)$. Donc ${\cal R}(\alpha)_{\vline
D_a}={p_1}_*({s_a}_*(\alpha))=(p_1\circ s_a)_*(\alpha)={p_a}_*(\alpha).$
\CQFD

On choisit un syst\`eme de coordonn\'ees affines de $\P_N$ $(x_1,\dots,x_n,y), n=N-1$, tel que l'hyperplan \`a
l'infini coupe le support de $\alpha$ proprement.
On lui associe un syst\`eme de coordonn\'ees sur
$G(1,N)$ en \'ecrivant les \'equations des droites sous la forme :
$x_i=a_iy+b_i$. Les droites pouvant s'\'ecrire de cette mani\`ere sont les droites ne rencontrant pas le
sous-espace de codimension deux $Y=Z=0$, o\`u $(X_1,\dots,X_n,Y,Z)$ est un syst\`eme de coordonn\'ees projectives
associ\'e.
On pose $x:=(x_1,\dots,x_n)$, $a:=(a_1,\dots,a_n)$.

En coordonn\'ees affines, on a $p_a(x,y):=x-ay$. 

\begin{lem} 
${\cal R}(\alpha)=\sum_I{u_{card(I)}(a,b)da^I\wedge db^{I^c}}$,
 avec: $$u_k(a,b):=\sum_{P_i}{Res_{P_i}{\alpha y^k/(l_1\dots l_n)}},$$ et les $P_i$ sont les diff\'erents points
d'intersection de la droite $\Delta_{a,b}$ avec le support de $\alpha$, et $l_i(x,y):=a_iy+b_i-x_i$. 
\end{lem}

\dem Ecrivons 
 que ${\cal R}(\alpha)$ est l'image directe par $p_1: U^*\times U\to U^*$ du courant l. r. $res_{I_U}{\alpha\wedge dl_1/l_1\wedge\dots\wedge dl_n/l_n}$,
 avec  $l_i(x,y)=a_iy+b_i-x_i$. D'apr\`es ce qu'on a vu du calcul de l'image 
directe, on obtient l'expression annonc\'ee.
\CQFD

{\it D\'emonstration du th\'eor\`eme $3$ : Premi\`ere \'etape.}

On suppose momentan\'ement que :

i) ${\cal R}(\alpha)$ se prolonge en une forme {\it holomorphe} sur $U^*$, et 

ii) ${\tilde U}^*$ v\'erifie la condition topologique suivante : 
les sections ${\tilde D}_a$ sont contractibles.

\begin{lem} $u_k(a,b)$ se prolonge holomorphiquement dans 
${\tilde D}_a$ pour tout $k$. 
\end{lem} 

\dem 
On a sur $U^*$ :
 ${\cal R}(\alpha y^k)=\sum_{I\subset\{ 1,\dots,n\} }{u_{k+ card(I)}da^I\wedge db^{I^c}}$, pour tout $k\ge 0$. 
 Comme ${\cal R}(\alpha y^k)$ est $d-$ferm\'ee sur $U^*$, on a ${\partial\over{\partial 
b_i}}u_{k+n}={\partial\over{\partial a_i}}u_{k+n-1}$, pour tout $i, 1\le i\le n$. 

L'hypoth\`ese que ${\cal R}(\alpha)$ se prolonge holomorphiquement nous donne d\'ej\`a que $u_0,\dots, u_n$ se
prolongent holomorphiquement dans
${\tilde  U}^*$; donc en particulier, sur ${\tilde D}_a$ pour tout $a$. Supposons que l'on ait prolong\'e $u_{n+k}$ 
sur ${\tilde D}_a$, pour un certain entier $k\ge 0$.

Consid\'erons la forme diff\'erentielle holomorphe sur ${\tilde
D}_a$: $$\phi_k:=\sum_{i=1}^n{\partial_{a_i}{u_{n+k}}db_i}.$$ Elle est $d-$ferm\'ee sur $D_a$, donc sur ${\tilde D}_a$.
D'apr\`es l'hypoth\`ese de contractibilit\'e,
 elle est donc $d-$exacte : $\phi_k=dv_k$, avec $v_k$ une fonction holomorphe sur ${\tilde D}_a$. Comme sur $D_a$,
elle est
\'egale
\`a
$du_{n+k+1}$, on en d\'eduit que
$u_{n+k+1}$, qui est holomorphe sur $D_a$, se prolonge holomorphiquement sur ${\tilde D}_a$ ($v_k+Cst$).
Par r\'ecurrence, tous les $u_k(a,b), k\ge 0$ se prolongent holomorphiquement sur ${\tilde D}_a$. 

\CQFD 

\begin{lem} Sous les hypoth\`eses pr\'ec\'edentes, $\alpha$ se prolonge \`a $\tilde U$ en courant l.
r. .
\end{lem}

\dem
 D'apr\`es ce qui
pr\'ec\`ede, les $u_k(a,b)={p_a}_*(\alpha y^k)$ se prolongent holomorphiquement dans les  domaines 
${\tilde D}_a$. D'apr\`es le th\'eor\`eme sur les traces, la 
restriction de $\alpha$ \`a $U_a$ se prolonge \`a 
${\tilde U}_a$, en un courant
$\alpha_a$. 

D'apr\`es le lemme ci-dessous, on a ${\alpha_a}_{\vline U\cap {\tilde U}_a}=\alpha_{\vline  U\cap {\tilde U}_a}$, donc on
peut en fait consid\'erer le courant $\alpha_a$ sur l'ouvert lin\'eairement $1-$concave $U\cup {\tilde U}_a$.

Consid\'erons la famille des couples $(V,\alpha_V)$, avec $V\subset {\tilde U}$ un ouvert lin\'eairement
$1-$concave contenant
$U$, et $\alpha_V$ un prolongement l.r. de $\alpha$ sur $V$, avec la relation d'ordre naturelle (si $V\subset V'$,
alors la restriction de $\alpha_{V'}$ \`a $V$ co\"\i ncide avec $\alpha_V$). Alors, cette famille admet un
\'el\'ement maximal, d'apr\`es le lemme de Zorn.

Soit $(W,\alpha_W)$ un tel \'el\'ement maximal. D'apr\`es ce qui pr\'ec\`ede, tout c\^one ${\tilde U}_x$ \`a
sommet $x$ dans $W$ est contenu dans $W$.

Consid\'erons un ouvert lin\'eairement $1-$concave $W\subset \tilde U$, 
tel que pour tout $x\in W$, le c\^one ${\tilde U}_x$ soit inclus dans $W$. Alors $V=\tilde U$. En effet, $W^*$ est
ouvert dans ${\tilde U}^*$. Mais $W^*$ est de plus ferm\'e dans ${\tilde U}^*$. Soit en effet un point $t\in
{\tilde U}^*$ sur la fronti\`ere de $W^*$, correspondant \`a une droite $\Delta_t$. Alors, $\Delta_t$ doit
rencontrer $W$ (sinon on pourrait trouver un voisinage ouvert de $t$ disjoint de $W^*$).
Soit $P\in \Delta_t\cap W$. Alors, $W$ contient le c\^one form\'e des droites passant par $P$ et contenues dans
$\tilde U$, et donc en particulier $\Delta_t$. Par cons\'equent, $t\in W^*$.
Comme ${\tilde U}^*$ est connexe, on a donc $W^*=U^*$, donc $W=\cup_{t\in W^*}=\tilde U$.
\CQFD 

\begin{lem} ${\alpha_a}_{\vline U\cap {\tilde U}_a}=\alpha_{\vline  U\cap {\tilde U}_a}$.
\end{lem}
\dem
Il s'agit de montrer que le courant $\alpha_a$, d\'efinit sur ${\tilde U}_a$, ne d\'epend pas de $a$ :
$\alpha_a=\alpha_{a'}$ sur ${\tilde U}_a\cap {\tilde U}_{a'}$.
Rappelons que $\alpha_a$ est le r\'esidu de la forme m\'eromorphe $\Psi_a:=r_a(x,y)dx\wedge dy/P_a(x,y)$.
Si on consid\'ere les $u_k$ et les $a_i$ qui entrent dans l'expression de $\Psi_a$ comme des fonctions m\'eromorphes sur
${\tilde U}^*$, on voit qu'on peut consid\'erer $\Psi_a$ comme une forme m\'eromorphe sur $I_{\tilde U}$.
Comme par ailleurs $\overline\partial[\Psi_a]=p_2^*(\alpha)$ sur $I_U$, on en d\'eduit qu'on a aussi
$\overline\partial[\Psi_a]=p_2^*(\tilde\alpha)$ sur $I_{\tilde U}$, pour un courant l.r. $\tilde\alpha$ sur $\tilde U$. 
Cela signifie que lorsque l'on change $a$, le r\'esidu de la forme m\'eromorphe
$\Psi_a$, $\overline\partial[\Psi_a]=\alpha_a$, ne change pas. 
\CQFD

{\it Deuxi\`eme \'etape : supression de la condition topologique sur $\tilde U$.}

On fait le m\^eme raisonnement que ci-dessus. Supposons donn\'e un ouvert $V\subset {\tilde
U}^*$ tel qu'on ait prolong\'e le courant $\alpha$ comme courant l.r. sur $\cup_{t\in
V}{\Delta_t}$. Alors si $V\not={\tilde U}^*$, on prend un point $t$ de la fronti\`ere de $V$
dans
${\tilde U}^*$, un voisinage ouvert $V_t$ de $t$ dans ${\tilde U}^*$, dont les sections sont
contractibles, et on prolonge
$\alpha$ comme courant l.r. dans $V_t':=\cup_{s\in V_t}{\Delta_t}$. On a donc prolong\'e
$\alpha$ sur $(V\cup V_t)'$. D'apr\`es Zorn, il existe un ouvert $V$ maximal, i.e. tel qu'on ne puisse plus
prolonger $\alpha$ sur $W'$, avec $V\subset W\subset {\tilde U}$. Un tel ouvert est n\'ecessairement
\'egal \`a ${\tilde U}$.

 {\it Troisi\`eme
\'etape : cas du prolongement m\'eromorphe.}

 On suppose maintenant juste que ${\cal R}(\alpha)$ se prolonge m\'eromorphiquement. Alors, l'ensemble
polaire $Z^*$ de ${\cal R}(\alpha)$ dans
${\tilde U}^*$ correspond aux droites rencontrant un certain sous-ensemble analytique $Z$ de
codimension deux dans $\tilde U$. ${\tilde U}-Z$ est lin\'eairement
$1-$concave,
${{\tilde U}-Z}^*={\tilde U}^*-Z^*$ est connexe, on peut donc appliquer le th\'eor\`eme dans le cas du prolongement
holomorphe, et trouver un prolongement $\tilde\alpha$ l.r. de $\alpha$ dans ${\tilde U}-Z$.
Le th\'eor\`eme de Remmert-Stein permet de prolonger le support $\tilde Y$ de $\tilde\alpha$ \`a travers $Z$; on
notera le prolongement encore $\tilde Z$.  Il reste donc \`a prolonger $\tilde\alpha$ \`a travers $\tilde Z$,
gr\^ace \`a notre hypoth\`ese de prolongement m\'eromorphe. 

\begin{lem} Soit $D$ un domaine de $\C^n$. Consid\'erons une fonction holomorphe $f$ en dehors d'une hypersurface
analytique $T\subset D$. Si $df$ est m\'eromorphe sur $D$, $f$ se prolonge m\'eromorphiquement sur $D$. 
\end{lem}
\dem 
La question est locale. D'apr\`es le th\'eor\`eme de Hartogs, il suffit de prolonger m\'ero\-morphiquement $f$ au voisinage
des points r\'eguliers de $T$. On peut donc supposer $T$ de la forme $\{y=0\}$, avec des coordon\'ees affines
$(x_1,\dots,x_{n-1},y)$. On peut supposer $f$ holomorphe dans une couronne $K:=D\times\{\vline
z\vline <\epsilon\}-\{\vline z\vline \le\eta,\eta<\epsilon\}$. Alors, $f$ peut s'\'ecrire sur $K$ comme
$f=f_+-f_-$, avec $f_+=\sum_{k\ge 1}{u_{-k}(x)y^{k-1}}, f_-=\sum_{k\ge 0}{u_k(x)/y^{k+1}}$,
avec $$u_k(x):=1/2i\pi\int_{\vline y\vline=\rho}{f(x,y)y^{k}} (\eta<\rho<\epsilon)$$ holomorphes sur $D$. Dire que $f$
se prolonge m\'eromorphiquement sur
$D\times\{\vline z\vline\le\epsilon\}$ revient \`a dire que $y^m f$ se prolonge holomorphiquement pour un certain
$m$, donc que $f_-$ comporte un nombre fini de termes non identiquement nuls. 
Remarquons gr\^ace \`a la convergence uniforme, on peut d\'eriver sous le signe somme. 
On obtient donc $f_y(x,y)=\sum_{k\in \Z}{u_{-k-1}(x)k y^{k-1}}$ sur $K$. 

Supposons que $df$ se prolonge m\'eromorphiquement sur $D\times\{\vline z\vline\le\epsilon\}$. Alors, il en est de
m\^eme de $f_z$. On en d\'eduit que seuls un nombre fini des $u_i,i\ge 0$ sont non identiquement nuls. Par
cons\'equent, $f$ aussi se prolonge m\'eromorphiquement.
\CQFD

Le lemme pr\'ec\'edent nous permet de montrer, suivant la technique employ\'ee ci-dessus, que les $u_k(a,b)$ se
prolongent m\'eromorphiquement sur $D_a$, pour tout entier $k\ge 0$.
Le lemme \ref{mero} permet alors de prolonger $\tilde\alpha$ \`a travers $Z$ dans $\tilde U$, en utilisant
plusieurs centres de projection. Ceci termine la d\'emonstration.
\CQFD

\begin{cor} Supposons ${\cal R}(\alpha)$ rationnel. Alors, $\alpha$ se prolonge en un courant l.r. sur $\P_N$.
\end{cor}

\section{Une variante d'un th\'eor\`eme de Rothstein.} 

Consid\'erons deux boules ouvertes dans l'espace $\C^N$, 
 contenues l'unes dans l'autre : $B\subset B'$. On se donne un sous-ensemble analytique 
de dimension pure $n\ge 2$ dans $B'-B$. Le th\'eor\`eme de Rothstein dit que ce 
sous-ensemble 
se prolonge en un sous-ensemble analytique de $B'$. 

Montrons comment le th\'eor\`eme pr\'ec\'edent d'inversion 
 pour la transformation d'Abel-Radon permet de montrer une variante du th\'eor\`eme 
de Rothstein. 

Cette variante est une g\'en\'eralisation du lemme \ref{pseudoconcave} 
ci-dessus. 

\begin{thm} 
Soit $U$ un domaine lin\'eairement $1-$concave de l'espace projectif $\P_N$. Supposons que $U$ contient un $2-$plan
complexe. On se donne un courant 
$\alpha$ l. r. de type 
$(N,1)$ dans $U$. Alors, $\alpha$ se prolonge dans $\P_N$ en courant l.r. . 
\end{thm} 

Nous aurons besoin de la notion de pseudoconcavit\'e au sens d'Andreotti.
 Un vari\'et\'e
analytique connexe $X$ est pseudoconcave au sens d'Andreotti s'il existe un domaine $U\subset X$ tel que pour tout
point de la fronti\`ere $x\in\partial U$, pour tout voisinage ouvert $U_x$ de $x$ dans $X$, toute fonction
holomorphe sur
$U_x-{\overline U}$ est holomorphe en $x$.
Il suffit pour v\'erifier cette propri\'et\'e, d'apr\`es Hartogs, d'exhiber un disque analytique $D$ de $X$, tel que
$D\cap
\overline{U}=\{x\}$. 

\begin{lem} Soit $U\subset \P_N$ un domaine lin\'eairement $1-$concave. Supposons que $U$ contienne un $2-$plan
complexe. Alors si l'ouvert $U^*\subset G(1,N)$ est connexe, il est pseudoconcave au sens d'Andreotti.
\end{lem} 

\dem On consid\`ere les $U_\epsilon$ d\'efinis 
par $\vline x_0\vline^2+\dots+\vline x_{N-3}\vline^2-\epsilon(\vline 
x_{N-2}\vline^2+\vline 
x_{N-1}\vline^2+\vline 
x_{N}\vline^2)<0$, o\`u $x_0=x_1=\dots=x_{N-3}=0$ sont les 
\'equations du $2-$plan. Ces domaines forment pour $\epsilon>0$ un syst\`eme fondamental de voisinages 
du $2-$plan. Soit $U_\eta$ un tel domaine contenu dans $U$.

Consid\'erons le domaine $V:=U_\eta^*\subset G(1,N)$. On se donne un point de la fronti\`ere $t\in\partial V$.
Alors, la droite correspondante n'est pas contenue dans $U_\eta$. Elle coupe donc la fronti\`ere, en un certain
point $x\in \overline{U_\eta}$. On consid\`ere un $2-$plan $H$, contenant la droite $\Delta$, et contenu dans
$\overline{U_\eta}$. On se donne un point $y$ de $H$ en dehors de $\Delta$. Alors, on consid\`ere les droites $D_z$
reliant $y$ et $z$, pour $z\in \Delta$. Alors, pour $z\not= x$, $D_z$ est contenu dans $U_\eta$, puisque
$H-\{x\}\subset U_\eta$. On a donc un $\P_1\subset \overline{V}$, qui ne rencontre $V$ qu'au point $t$, ce qui
montre la pseudoconcavit\'e, d'apr\`es ce qu'on a vu. 
\CQFD

{\it D\'emonstration du th\'eor\`eme.} 
Soit $U$ un domaine de l'espace projectif $\P_N$, contenant un $2-$plan complexe.
La transform\'ee ${\cal R}(\alpha)$ est m\'eromorphe sur $U^*$. 
P. Dingoyan (\cite{Dingoyan}) a montr\'e que toute fonction m\'eromorphe sur un ouvert pseudoconcave d'un
vari\'et\'e alg\'ebrique est rationnelle.
 Comme $U^*$ est pseudoconcave au sens d'Andreotti, on en d\'eduit que ${\cal R}(\alpha)$ est rationnelle.
 D'apr\`es le corollaire ci-dessus, $\alpha$ se prolonge sur 
$\P_N$ en un courant localement r\'esiduel. 
\CQFD 

\section{Questions ouvertes.}

1. On pourrait g\'en\'eraliser le th\'eor\`eme $3$, en ne supposant pas a priori que ${\cal R}(\alpha)$ se prolonge dans un
domaine de la forme ${\tilde U}^*$, pour $\tilde U\subset \P_N$ contenant $U$. Le th\'eor\`eme deviendrait alors : si
${\cal R}(\alpha)$ se prolonge dans un doamine $D$ contenant $U^*$, alors $\alpha$ se prolonge dans $D':=\cup{t\in
D}{\Delta_t}$ comme courant l.r..
 A fortiori, ${\cal R}(\alpha)$ se prolongerait alors dans l'"enveloppe" de $D$, i.e.
${D'}^*$.

2. Si $B$ est juste pseudoconvexe dans $\P_3$ (mais pas $\C-$convexe), 
 on aimerait 
avoir un exemple tel que 
$\alpha$ de bidegr\'e $(3,1)$ dans le compl\'ementaire de $B$ qui ne se prolonge pas. 
$B$ doit rencontrer toutes hypersurfaces alg\'ebriques. Il suffirait de trouver une 
hypersurface de $U=\P_3-B$ non alg\'ebrique. 

3. On se donne un ouvert lin\'eairement $2-$concave $U\subset \P_N$ un courant l.r. de bidegr\'e $(N,2)$ dans $U$.
On suppose ${\cal R}(\alpha)=0$. On pourrait d\'eduire de ce qui pr\'ec\`ede l'existence d'un prolongement l.r.
$\tilde\alpha$ \`a $\P_N$, si l'on pouvait montrer que dans $U$, le fait que $\alpha$ est $\overline\partial-$exact
implique $\alpha=\overline\partial\beta$, avec $\beta$ l.r. dans $U$.

\hspace{60pt} {22, rue Emile Dubois, 75014 PARIS, France}
 
\hspace{60pt} {bruno.fabre@iecn.u-nancy.fr}

\end{document}